\documentclass[11pt]{article}
\vspace{7cm} \textwidth 14cm \textheight 21.6cm
\usepackage{amsmath}
\usepackage{amsfonts}
\usepackage{latexsym,amsthm,amscd}
\newtheorem{thm}{Theorem}
\newtheorem{cor}{Corollary}
\newtheorem{prop}{Proposition}

\newtheorem{lem}{Lemma}
\newtheorem{Def}{Definition}

\newcounter{alphthm}
\newtheorem{lemm}[alphthm]{Lemma}

\newcommand{\be}{\begin{equation}}
\newcommand{\ee}{\end{equation}}
\newcommand{\bc}{\begin{cor}}
\newcommand{\ec}{\end{cor}}
\newcommand{\ben}{\begin{enumerate}}
\newcommand{\een}{\end{enumerate}}
\newcommand{\beq}{\begin{eqnarray}}
\newcommand{\eeq}{\end{eqnarray}}
\newcommand{\beqn}{\begin{eqnarray*}}
\newcommand{\eeqn}{\end{eqnarray*}}

\newcommand{\bpf}{\begin{proof}}
\newcommand{\epf}{\end{proof}}
\newcommand{\bl}{\begin{lem}}
\newcommand{\el}{\end{lem}}
\newcommand{\bp}{\begin{prop}}
\newcommand{\ep}{\end{prop}}
\newcommand{\bd}{\begin{Def}}
\newcommand{\ed}{\end{Def}}
\newcommand{\bt}{\begin{thm}}
\newcommand{\et}{\end{thm}}
\newcommand{\R}{I\!\! R} 

\title{Strong rigidity of  constant curvature Finsler manifolds \thanks{
   \hspace{4cm} \emph{Keywords and phrases }: Finsler, rigidity,  constant curvature, second order differential equation, adapted coordinates. }}

\author{\small Asanjarani A.  Bidabad B.
\\ \small Amirkabir University of Technology (Tehran
Polytechnic), \\ \small  Faculty of Mathematics and Computer Sciences,\\ \small 424
Hafez Ave., 15914 Tehran, Iran.\\ \small  bidabad@aut.ac.ir }
\date{}

\begin{document}
\maketitle
\begin{abstract}
Here, an extension of the Obata-Tanno's theorem to Finsler geometry is established and the following rigidity result is obtained; Every complete connected Finsler manifold  of positive constant flag curvature is isometrically homeomorphic  to an $n$-sphere
 equipped with a certain Finsler metric,  and vise versa.

\end{abstract}
%
\section{Introduction.}
\label{intro}Rigidity describes  quite different concepts in
 geometry. Historically, the first rigidity
theorem, proved by Cauchy in 1813, states that if the faces of a
convex polyhedron were made of metal plates and the edges were
replaced by hinges, the polyhedron would be rigid \cite{Cau}.
Although rigidity problems were of immense interest to engineers,
intensive mathematical study of these types of problems has occurred
only relatively recently \cite{Con}. Geometrically, sometimes an
object is said to be rigid if it has flexibility and not elasticity,
that is to say, invariant by each isometry. In Riemannian geometry
the sectional curvature is invariant under isometries. Hence a space
of positive constant curvature is transformed into the same space by
each isometry. This fact is sometimes described as the "strong
rigidity" of a space of constant curvature. So far, in  Finsler
geometry, the encountered rigidity results are rather slightly
weaker and they say that under such and such assumptions about the
flag curvature -analogous to the sectional curvature in Riemannian
geometry- the underlying Finsler structure must be either Riemannian
or locally Minkowskian. A famous treatise of this kind due to
Akbar-Zadeh \cite{AZ} in 1988 who has established in compact case
the following rigidity theorem; \emph{ Let $(M,g)$  be a compact
without boundary Finsler manifold of constant flag curvature $K$.
 If $K<0$, then $(M,g)$ is Riemannian.
 If $K=0$, then $(M,g)$ is locally Minkowskian.}

 In 1997 Foulon \cite{Fo1} addresses the case of strictly negative flag curvature by imposing
the additional hypothesis that the curvature be covariantly constant
along a distinguished vector field on the homogeneous bundle of
tangent half lines. He shows that, under these conditions, as in the
Akbar-Zadeh's theorem, the Finsler structure is Riemannian. Next in
2002, Foulon came back to this problem and  presented a strong
rigidity theorem for symmetric compact Finsler manifolds with
negative curvature and proves that it is isometric to a locally
symmetric negatively curved Riemannian space \cite{Fo2}. This
extends the Akbar-Zadeh's rigidity theorem, to a so called, strong
rigidity one.

 Shen  in 2005 \cite{Sh1} observes the case of negative but not necessarily
constant flag curvature by imposing the additional hypothesis that
the $S$-curvature be constant and proved that the Akbar-Zadeh's
rigidity theorem still holds.

In 2007 following the several rigidity theorems in the two joint
papers \cite{KY1} and \cite{KY2},  Kim  proved that \cite{Ki}
\emph{; any compact locally symmetric Finsler manifold with positive
constant flag curvature is Riemannian. }

One of the present authors has also established in   a joint paper in 2007, some rigidity theorems as an application of connection theory in Finsler geometry \cite{BT}.

Here, we study the strong rigidity of Finsler manifolds of positive
constant  curvature  and show that;
\emph{Let $(M,g)$
be an $n$-dimensional  complete connected Finsler manifold. Then it is  of positive constant flag curvature
   $K>0$ if and only if  $(M,g)$ is isometrically homeomorphic  to an $n$-sphere
 equipped with a certain Finsler metric.}
Particularly,  this result completes the
Akbar-Zadeh's rigidity theorem for positive constant flag curvature.
Meanwhile, we have extended the Obata-Tanno's theorem to Finsler
geometry as follows.
\emph{Let   $(M,g)$  be a complete connected Finsler manifold of
  dimension $n\geq2$.  In
  order that  there is  a non-trivial solution
 of
 $
 \nabla^{H}\nabla^{H}\rho+ C^2 \rho g=0
 $
  on $M$, it is necessary and
sufficient that $(M,g)$ be isometric to an n-sphere of radius $1/C$.}

This result  leads us  to illustrate a definition for an $n$-sphere in
Finsler geometry.

\section{Preliminaries.}
 Let $(M,g)$ be an n-dimensional Finsler
manifold and  $TM\rightarrow M$ the  bundle
  of its tangent vectors. Using the local coordinates $(x^{i},y^{i})$ on $TM$, called the line
elements, we have the local field of frames
$\{\frac{\partial}{\partial x_{i}},\frac{\partial}{\partial
y_{i}}\}$ on $TTM$.

 Let $\rho:M\rightarrow \R$  be a scalar
function on $M$ and consider the following second order differential
equation
 \be
 \label{c-field}
\nabla^{H}\nabla^{H}\rho= \phi g,
 \ee
where $\nabla^{H}$ is the Cartan horizontal covariant derivative and
$\phi$ is a function of $x$ alone, then we say that Eq. (\ref{c-field}) has a solution $\rho$. The connected component of a
regular hypersurface defined by $\rho=constant$, is called a\emph{
level set of $\rho$}. We denote by $\verb"grad" \rho$ the gradient
vector field of $ \rho$ which is locally
 written in the form $\verb"grad" \rho =  \rho^{i}\frac{\partial}{\partial
x^i}$, where  $\rho^i = g^{ij} \rho_j$, $\rho_j =\frac{\partial
\rho}{\partial x^j}$ and $i, j, ...$ run over the range $1,...,n$.

The partial derivatives $\rho_j $ are defined on the manifold $M$
while $ \rho^{i}$ the components of $\verb"grad" \rho$ may be defined
on its slit tangent bundle $TM_0$. Hence, $\verb"grad" \rho$
 can be considered as a section of $\pi^*TM\rightarrow TM_0$,
  the pulled-back tangent bundle over $TM_0$,
and  its trajectories lie on $TM_0$.

One can easily verify that the canonical projection of the trajectories of the vector field
$\verb"grad" \rho $ are geodesic arcs on $M$ \cite{AB}. Therefore, we can choose a
local coordinates
  $(u^1=t, u^2,...,u^n)$ on $M$  such that $t$ is the parameter of the geodesic containing
  the  projection of a trajectory of the
vector field $\verb"grad" \rho$ and the level sets of $\rho$ are given by $t=$constant. These geodesics are called \emph{
$t$-geodesics}. Since in this local coordinates,
the level sets of $\rho$ are given by $t=$constant, then $\rho$ may be considered as a function of $t$ only.
 In the sequel we will refer to these level sets and this local coordinates as \emph{$t$-levels} and \emph{adapted
coordinates}, respectively.

 Let $(M,g)$ be a Finsler manifold and $\rho$  a
non-trivial solution of Eq. (\ref{c-field})  on $M$. Then in an adapted coordinates, components of
  the Finsler metric tensor  $g$ are given by
 $$\footnotesize
(g_{ij})= \left(\begin{array}{lcl}
1 \quad {0\quad \ldots\qquad{ 0}}\\
0 \quad {g_{22} \quad\ldots\quad g_{2n}}\\
\vdots \qquad{ \quad \ldots}\\
0 \quad {g_{n2} \quad\ldots \quad g_{nn}}
\end{array}\right),
 $$
and $t$ may be regarded as
 the arc-length parameter of $t$-geodesics.
It can be easily verified that the Finsler metric form of $M$ is given by
  \be \label{meter}
ds^{2}=(dt)^{2}+ \rho'^{2}f_{\gamma\beta}du^{\gamma}du^{\beta},
  \ee
  where $f_{\gamma\beta}$ given by $g_{\gamma\beta}=\rho'^{2}f_{\gamma\beta}$ are components of a Finsler metric tensor
 on a $t$-level of $\rho$ and  $g_{\gamma\beta}$ is the induced
 metric tensor of this $t$-level. Here and every where in this
 paper, the Greek indices $\alpha, \beta, \gamma, ...$ run over the
 range $2,3, ..., n.$

If $g(\verb"grad" \rho,\verb"grad" \rho)=0$ in some points of $M$,
then $M$ possesses some interesting properties.
 A point $o$ of $(M,g)$ is called a \emph{critical point} of $\rho$ if the
vector field $\verb"grad" \rho$ vanishes at  $o$, or equivalently if $
\rho'(o)=0$.

 \setcounter{alphthm}{0}
  \begin{lemm}\label{lem} \cite{AB}
 Let $(M,g)$ be  an n-dimensional Finsler manifold which admits a non-trivial
solution
 $\rho$  of Eq. (\ref{c-field}). If $\rho$ admits a critical
point, then any $t$-level $\overline{M}$ with
Finsler metric form $\overline{ds}^{2}= f_{\gamma\beta}du^{\gamma}
du^{\beta}$ has the positive
constant flag curvature  $\rho''^2(0)$.
 \end{lemm}

\section{A special solution.}
Let $(M,g)$ be an n-dimensional Finsler manifold and
$\rho:M\rightarrow \R$  a solution of Eq. (\ref{c-field}). If  $\phi$ is a linear  function of $\rho$
with constant coefficients,  then we say that $\rho$ is a
\emph{special solution} of Eq. (\ref{c-field}). Hence, any
special solution of Eq. (\ref{c-field}) can be written in the form
 \be
 \label{sc-field}
\nabla^{H}\nabla^{H}\rho= (- K \rho + B) g,
 \ee
where $K$ and $B$ are constants.
 Along any geodesic with arc-length $t$,  Eq. (\ref{sc-field}) reduces to the
 ordinary differential equation
 \be \label{equation}
\frac{d^{2}\rho}{dt^{2}}=-K \rho + B .
 \ee
Now for the special case $K=C^2>0$ and $B=0$, we have
 \be
\label{sc-field1} \frac{d^{2}\rho}{dt^{2}}+ C^2\rho=0.
 \ee
 By a
suitable choice of the arc-length $t$, a solution of Eq. (\ref{sc-field1}) is given
by
 \be \label{rho1}
\rho (t)= A \cos Ct,
 \ee and its first derivative is
 \be
\rho' (t)=-AC \sin Ct.
 \ee
We can see at a glance, that it might appear   two critical points
corresponding to $t=0$ and $t=\frac{\pi}{C}$ on $M$, where these points are periodically repeated.
Hence, if $\rho$ is a non-trivial solution of
Eq. (\ref{sc-field1}), then it can be written in the following form
 \be \label{rho}
 \rho(t)=\frac{-1}{C}\cos Ct\ \quad( A=\frac{-1}{C}).
 \ee
Taking into account Eq. (\ref{meter}),
 the  metric form of $M$ becomes
  \be\label{metric form}
ds^2= dt^2+(\sin Ct)^2\overline{ds}^2,
 \ee
where $\overline{ds}^2$ is the metric form of a $t$-level of $\rho$
given by $\overline{ds}^2=f_{\gamma\beta}du^{\gamma}du^{\beta}$.

\section{Finsler manifolds of positive constant flag curvature.} Let
 $(x,y)$ be the line element of $TM$ and $P(y,X)\subset T_{x}(M)$  a 2-plane
  generated by the vectors $y$ and $X$ in
$T_{x}(M) $. Then the \emph{flag curvature} $K(x,y,X)$ with respect
to  the plane $P(y,X)$ at a point $x\in M$ is defined by
$$K(x,y,X):=\frac{g(R(X,y)y,X)}{g(X,X)g(y,y)-g(X,y)^{2}},$$
where $R(X,y)y$ is the $h$-curvature tensor of Cartan connection. If
$K$ is independent of $X$, then $(M,g)$ is called \emph{space of}
\emph{scalar curvature}. If $K$ has no dependence on $x$ or $y$,
then the Finsler manifold is said to be of \emph{constant curvature},
see for example \cite{AZ2}.  It can be easily verified that the
components of the $h$-curvature tensor of Cartan connection in
adapted coordinates are given by (see page 7 in \cite{AB})
 \be \label{2.21}\begin{array}{ccc}
R^\alpha_{\ 1 \gamma 1}= - R^\alpha_{\ \gamma 1
1}=(\frac{\rho'''}{\rho'})\delta_{\gamma}^{\alpha},\\\\
R^1_{\ 1 \gamma \beta}= - R^ 1_{\ \gamma 1 \beta}=-\rho' \rho'''
f_{\gamma\beta},\\\\
R^\alpha_{\ \delta \gamma \beta}=\overline{R}^\alpha_{\ \delta
\gamma \beta}-(\rho'')^2(f_{\gamma\beta}\delta_{\delta}^{\alpha}-
f_{\delta\beta}\delta_{\gamma}^{\alpha}),
 \end{array}
 \ee
where $\overline{R}^\alpha_{\ \delta
\gamma \beta}$ are components of $h$-curvature tensor related to the metric form $\overline{ds}^2$ on a $t$-level of $\rho$.
 \bp \label{constant curvature}
Let $(M,g)$ be an n-dimensional Finsler manifold. Then  it is of
constant flag curvature $K=C^2>0$, if and only if,
there is a  non-trivial solution of $\nabla^{H}\nabla^{H}\rho= (- C^2 \rho + B) g$ on $M$.
 \ep
 \bpf
 Let $(M,g)$ be a Finsler manifold, then it is of constant flag
curvature $K$, if and only if the components of the $h$-curvature
tensor are given by
 \be \label{K}
R^i_{\ hjk}= K (\delta^i_hg_{jk}-\delta^i_jg_{hk}).
 \ee
Using  Eq. (\ref{K}), one can easily drive the differential equation
 \be \label{bcs}
\ddot{A}+KAg=0,
 \ee
where $A$ is the Cartan torsion tensor, $\dot{A}_{ijk}:=(\nabla^H_sA_{ijk})y^s$ and
$\ddot{A}_{ijk}:=(\nabla^H_s\nabla^H_tA_{ijk})y^sy^t$ (see for example \cite{AZ2} or \cite{BCS}).

Fix any $X,
Y, Z \in \pi^*TM$ at $v \in I_xM=\{w \in T_xM,F(w)=1\}$. Let
 $c:I\!\!R
\rightarrow M$ be the unit speed geodesic in $(M,F)$ with
${{dc} \over {dt}}(0)=v$ and  $\hat c:={dc \over{dt}}$ be the
canonical lift of $c$ to $TM_0$. Let $X(t)$, $Y(t)$ and $Z(t)$
denote the parallel sections along $\hat c$ with $X(0)=X$, $Y(0)=Y$
and $Z(0)=Z$. Put $A(t)=A(X(t),Y(t),Z(t))$, $\dot A(t)=\dot
A(X(t),Y(t),Z(t))$ and $\ddot A(t)=\ddot A(X(t),Y(t),Z(t))$. Indeed
along geodesics, we have $ \frac{d  A}{dt}=\dot A$, $ \frac{d \dot A}{dt}=\ddot A$ and  Eq. (\ref{bcs}) becomes
\be \label{bcs2}
\frac{d^2A(t)}{dt^2}+K A(t)=0.
 \ee
The general solution of this differential equation is $A(t)=A_0 \cos \sqrt{K}t+B_0 \sin \sqrt{K}t$.
Therefore, Eq. (\ref{sc-field1}) which represents a special case of Eq. (\ref{sc-field}) along  geodesics, has a non-trivial solution on $M$.


Conversely, let $\rho$ given by Eq. (\ref{rho}) be a solution of Eq. (\ref{sc-field})  on $M$. Then
there is an adapted coordinate on $M$ for which the components of $h$-curvature are given by three equations (\ref{2.21}). Hence, first and second  equations  of
Eqs. (\ref{2.21}) satisfy
 \be\label{third }
 R^i_{\ hjk}=
\frac{-\rho'''}{\rho'}(\delta^i_hg_{jk}-\delta^i_jg_{hk}).
 \ee
   Differentiate  Eq. (\ref{rho}) and replace  the first and third derivatives
  of $\rho$, we obtain $\frac{-\rho'''}{\rho'}=C^2$. Therefore, Eq. (\ref{K}) is satisfied by two first equations of  Eqs. (\ref{2.21}).

 For checking the third equation of Eqs. (\ref{2.21}), we recall that as we saw in section 3,  $\rho$ has critical points on $M$. Thus,
subject to Lemma \ref{lem}, the
$t$-levels of $\rho$ are spaces of the positive constant curvature
$\rho''^2(0)=C^2$. Therefore,  the third equation of Eqs. (\ref{2.21}) becomes
$$
R^\alpha_{\ \delta \gamma \beta}=(C^2 -\rho''^2)(f_{\gamma\beta}\delta_{\delta}^{\alpha}-
f_{\delta\beta}\delta_{\gamma}^{\alpha}).
$$
Substituting  $g_{\alpha\beta}= \rho'^2 f_{\alpha\beta}$ and first and second derivatives of  $\rho$,  we obtain
$$R^\alpha_{\ \delta \gamma \beta}=C^2(g_{\gamma\beta}\delta_{\delta}^{\alpha}-
g_{\delta\beta}\delta_{\gamma}^{\alpha}).$$
  Consequently, Eq. (\ref{K}) is satisfied by
all three components of Cartan $h$-curvature tensor in adapted
coordinates,   and the Finsler manifold $(M,g)$ is of constant flag curvature
$K=C^2$.
 \epf


Now, we are in a position to prove an extension of the \emph{Obata-Tanno's} theorem
for Finsler manifolds.
 \bt\label{isometer2}
 Let   $(M,g)$  be a complete connected Finsler manifold of
  dimension $n\geq2$.  In
  order that  there is  a non-trivial solution
 of
 \be \label{sc-field2}
 \nabla^{H}\nabla^{H}\rho+ C^2 \rho g=0
 \ee
  on $M$, it is necessary and
sufficient that $(M,g)$ be isometric to an n-sphere of radius $1/C$.
  \et
 \bpf
 If there is a non-trivial
 solution of  Eq. (\ref{sc-field1}) on $(M,g)$, then according to Proposition \ref{constant curvature}, it is of positive constant
 flag curvature $C^2$. Preceding argument in the section $3$ shows that the  metric form of
 $(M,g)$ is given by Eq. (\ref{metric form}), which is the polar form of a Finsler
metric on a standard sphere of radius $1/C$ \cite{Sh2}. Therefore, we have the sufficiency.

Conversely, if $(M,g)$ is isometric to an n-sphere of radius
$\frac{1}{C}$, then the metric form of  $M$ is given by
$ds^{2}=(dt)^{2}+ \sin ^{2}(Ct)\overline{ds}^2$, where
$\overline{ds}^2$ is the metric form of a $t$-level of $M$.
This is the  polar form of a
 Finsler
 metric on an $n$-sphere in $I\!\!R^{n+1}$ with the positive  constant curvature $C^2$ .
Now if we substitute  the derivative of $\rho(t)=-1/C \cos Ct$, in the metric form of $M$, we obtain  $ds^{2}=(dt)^{2}+ \rho'(t)\overline{ds}^2$, where
$\rho(t)$ is a non-trivial solution of the second order differential equation
(\ref{sc-field1}), or equivalently (\ref{sc-field2}) along geodesics.
 \epf
Now taking into account the number of critical points of $\rho$,
we have
\bc \label{iff}
Let $(M,g)$ be a complete connected Finsler manifold. In order that
there is a non-trivial solution of $\nabla^{H}\nabla^{H}\rho+ C^2\rho
g=0$, it is necessary and sufficient that $(M,g)$ be isometrically
homeomorphic to an $n$-sphere.
 \ec
 \bpf
Let $(M,g)$ admit a non-trivial solution of  $\nabla^{H}\nabla^{H}\rho+ C^2\rho
g=0$,
then from Theorem \ref{isometer2} we know that it is isometric to an n-sphere of
radius $1/C$.
On the other hand, $M$ is complete and  Proposition \ref{constant curvature} says that, the Finsler manifold $(M,g)$ is of positive constant curvature. Therefore, by means of  extension   of the Meyers's theorem  for Finsler manifolds \cite{AZ}, $M$ is compact.
Thus, the function $\rho$ admits its absolute maximum and
minimum values on $M$. Consequently, $\rho$ has two critical points on $M$ and from  extension of
the Milnor's theorem for Finsler manifolds \cite{L}, $(M,g)$ is homeomorphic
to an $n$-sphere.

 Conversely, let $(M,g)$ be isometrically homeomorphic to an $n$-sphere of radius $1/C>0$, then corresponding to
  Theorem \ref{isometer2}, there is a non-trivial solution
 of $\nabla^{H}\nabla^{H}\rho+ C^2\rho
g=0$, on  $M$.
 \epf
 Following the Obata-Tanno's theorem in Riemannian geometry a unit
sphere is characterized by existence of solution of the differential
equation $\nabla\nabla f + f g=0$, where $f$ is a certain
function on $M$ and $\nabla$ is the Levi-Civita connection associated to the Riemannian metric $g$ \cite{G}. Similarly, Theorem \ref{isometer2} shows
that in Finsler geometry a unit sphere is characterized by existence
of solution of  $\nabla^{H}\nabla^{H}\rho+ \rho g=0$, where $\rho$
is a certain function on $M$. In analogy with  Riemannian geometry,
this leads to a definition for an $n$-sphere in Finsler geometry as
follows;
  \emph{Let $(M,g)$ be a complete connected Finsler manifold of positive constant flag curvature,
then it is said to be a {\it Finslerian $n$-sphere}, if it is
isometrically homeomorphic to an $n$-sphere}. \footnote{ That is, $M$
is homeomorphic to an $n$-sphere and $g$ is isometric to a certain
Finsler metric on an $n$-sphere.}
\bt
Let $(M,g)$ be an $n$-dimensional complete connected Finsler manifold. Then it is  of positive constant flag
curvature $K=C^2$, if and only if, $(M,g)$ is isometrically homeomorphic to an
$n$-sphere of radius $1/C$ endowed with a certain Finsler metric.
 \et
 \bpf
  Let $(M,g)$ be of positive constant curvature $C^2$. As a consequence of
  Proposition \ref{constant curvature}, There is  a non-trivial solution of Eq. (\ref{sc-field1}) on $M$
  and by means of Corollary \ref{iff} it is isometrically homeomorphic to an
$n$-sphere of radius $1/C$ equipped with a certain Finsler metric form $ds^{2}=(dt)^{2}+ \sin ^{2}(Ct)\overline{ds}^2$.

 Conversely, let $(M,g)$ be a  Finsler manifold which is isometrically homeomorphic to an $n$-sphere of radius $1/C>0$. Then, by means of Corollary \ref{iff},  $M$
admits a non-trivial solution
 of Eq. (\ref{sc-field1}). It follows from Proposition \ref{constant
 curvature} that $(M,g)$ is of positive constant flag curvature $C^2$.
 \epf


\begin{thebibliography}{}
%
%

\bibitem{AZ2} Akbar-Zadeh H.,  Initiation to global Finsler
geometry, North Holland, Mathematical Library, Vol \textbf{68},
(2006).

\bibitem{AZ} Akbar-Zadeh H., Sur les spaces de Finsler \`{a} courbures
sectionnelles constantes,  Acad. Roy. Bull. Cl. Sci. (5)\textbf{ 74},
(1988), 281-322.


\bibitem{AB} Asanjarani A., Bidabad B., Classification of complete Finsler manifolds through a second order differential equation, Differential
 Geometry and its Application,to appear, (2007).

\bibitem{BCS} Bao D., Chern S.S., Shen Z., An introduction to Riemann-Finsler geometry, Springer-Verlag, (2000).
\bibitem{BT} Bidabad B., Tayebi A., A classification of some Finsler connections
and their applications, Publ. Math. Debrecen, \textbf{71}/3-4, (2007), 253-266.

\bibitem{Cau}Cauchy  A. L., XVIe Cahier IX,
\textbf{87E/89}, (1813).

\bibitem{Con} Connelly R.,  Ch. 1.7 in Handbook of Convex
Geometry, Vol. \textbf{A} (Ed. P. M. Gruber and J. M. Wills).
North-Holland, Amsterdam, Netherlands, (1993) pp. 223.

\bibitem{Fo1} Foulon P., Locally symmetric Finsler spaces in negative curvature, C. R. Acad. Sci. Paris S$\acute{e}$r. I
Math.,  no.\textbf{10}, 324 (1997) 1127-1132.

\bibitem{Fo2} Foulon P., Curvature and global rigidity in Finsler manifolds, Houston J. Math.,\textbf{ 28}, (2002), 263-292.


\bibitem{G} Gallot  S., \'{E}quations diff\'{e}rentielles
caract\'{e}ristiques de la sph\'{e}re, Ann. Scient. \'{E}c. Norm.
Sup. 4 s\`{e}rie, t. \textbf{12}, (1979), 235-267.


\bibitem{KY1} Kim  C. W. and Yim  J. W., Rigidity of noncompact Finsler manifolds, Geometriae Dedicata \textbf{81}, (2000),
 245-259.

\bibitem{KY2} Kim  C. W. and Yim  J. W., Finsler manifolds with positive constant flag curvature, Geometriae Dedicata \textbf{98}, (2003), 47-56.


\bibitem{Ki} Kim  C. W., Locally symmetric positively curved Finsler
spaces, Arch. Math. \textbf{88}, No.4, (2007), 378-384.



\bibitem{L} Lehmann  D.,   s\'{e}minaire Ehresmann,
Topologie et g\'{e}ometrie diff\'{e}rentielle,\textbf{ 6} (1964).



\bibitem{Sh1}
Shen Z., Finsler manifolds with nonpositive flag curvature and constant $S$-curvature, Math. Z. Volume no.\textbf{ 3}, 249 (2005),625-639.

\bibitem{Sh2} Shen Z., Differential Geometry of Spray and
Finsler Spaces, Kluwer Academic Publishers, Dordrecht (2001).


\end{thebibliography}
\end{document}